\documentclass[
  twoside,
  aps,
  pre]{revtex4}

\usepackage[usenames]{color}
\usepackage{amssymb,amsmath,amsthm}
\usepackage{graphicx}

\usepackage[colorlinks=true,
linkcolor=webgreen,
filecolor=webbrown,
citecolor=webgreen]{hyperref}

\definecolor{webgreen}{rgb}{0,.5,0}
\definecolor{webbrown}{rgb}{.6,0,0}

\RequirePackage[ngerman,american]{babel}

\usepackage{color}
\usepackage{float}

\newcommand{\seqnum}[1]{\href{http://www.research.att.com/cgi-bin/access.cgi/as/~njas/sequences/eisA.cgi?Anum=#1}{\underline{#1}}}

\newcommand{\E}{\mathrm{e}}
\newcommand{\bigo}{\mathcal{O}}

\begin{document}

\title{Asymptotics of Lagged Fibonacci Sequences}

\author{Stephan Mertens}
\email[URL:]{www.ovgu.de/mertens}
\affiliation{{\selectlanguage{ngerman}{Institut f"ur Theoretische Physik,
    Otto-von-Guericke Universit"at, PF 4120, 39016 Magdeburg, Germany}}\\
    Santa Fe Institute, 1399 Hyde Park Road, Santa Fe, NM 87501, USA}

\author{Stefan Boettcher}
\email[URL:]{www.physics.emory.edu/faculty/boettcher/}
\affiliation{Department of Physics, Emory University, Atlanta GA 30322, USA}

\begin{abstract}
  Consider ``lagged'' Fibonacci sequences $a(n) = a(n-1)+a(\lfloor
  n/k\rfloor)$ for $k > 1$.  We show that $\lim_{n\to\infty}
  a(kn)/a(n)\cdot\ln n/n = k\ln k$ and we demonstrate the slow
  numerical convergence to this limit and how to deal with this slow
  convergence. We also discuss the connection between two classical
  results of N.G. de Bruijn and K. Mahler on the asymptotics of
  $a(n)$.
\end{abstract}

\maketitle

\section{Introduction}

Let $k > 1$ be an integer and consider ``lagged'' Fibonacci type
sequences
\begin{equation}
  \label{eq:def-lagged-fibo}
  a_k(n) = a_k(n-1) + a_k\left(\left\lfloor\frac{n}{k}\right\rfloor\right)
\end{equation}
with initial value
\begin{equation}
  \label{eq:initial-values}
  a_k(0) = 1\,.
\end{equation}
These ``almost linear recurrence'' has many interesting arithmetical
properties \cite{knuth:66}.  The value $a_k(n)$ equals the number of
$k$-ary partitions of $kn$, and the corresponding sequences are listed
in the OEIS as \seqnum{A000123}, \seqnum{A005704}, \seqnum{A005705}
and \seqnum{A005706} for $k=2,3,4,5$, respectively.  In this
contribution we will study the asymptotical behavior of the ratio
\begin{equation}
  \label{eq:def-c}
  c_k(n) = \frac{a_k(kn)}{a_k(n)}\,\frac{\ln n}{n}\,.
\end{equation}

The OEIS entry for \seqnum{A000123} quotes a conjecture due to Benoit Cloitre, claiming that 
\begin{equation}
  \label{eq:cloitre}
  \lim_{n\to\infty} c_k(n) = \text{const.} = 1.63\ldots\,.
\end{equation}
The same conjecture (but with $\text{const.}=1.64\ldots$) appears for
the related sequence \seqnum{A033485}. We will prove that the
essential part of the conjecture (existence of the limit) is true, but
that its numerical part is incorrect. In particular, we will apply a
classical result of de Bruijn \cite{de-bruijn:48} to prove that
\begin{equation}
  \label{eq:our-c}
  \lim_{n\to\infty} c_k(n) = k\ln k\,.
\end{equation}
Note that $2\ln 2 = 1.386\ldots$, which differs significantly from
the value in \eqref{eq:cloitre}.

In the second part we will discuss the rate of convergence of
$c_k(n)$.  It turns out that this rate is so slow that
straightforward numerical measurements of $c_k(n)$ cannot be used for
an accurate measurement of $c_k(\infty)$. This may explain the inaccurate
numerical value in \eqref{eq:cloitre}. It turns out that another
classical result on the asymptotics of $a_k(n)$ due to K. Mahler
\cite{mahler:40} can be used as a device for an accurate numerical
determination of $c_k(n)$ all the way to the asymptotic regime.

In the final part we will discuss the connection between the two
asymptotic formulas of de Bruijn and Mahler.

\section{Asymptotics}

Using an integral representation (Mellin transformation) of the generating function for $a_k(n)$ and a saddle point integration, de Bruijn
\cite{de-bruijn:48} showed that
\begin{equation}
  \label{eq:de-bruijn-asymptotics}
  \begin{aligned}
  \ln a_k(n) = & \frac{1}{2\ln k}(\ln n -\ln\ln n)^2 
                +\left(\frac{1}{2}+\frac{1}{\ln k} + \frac{\ln\ln k}{\ln k}\right)\ln n 
               - \left(1+\frac{\ln\ln k}{\ln k}\right)\ln\ln n\\
               & + \left(1+\frac{\ln\ln k}{2\ln k}\right)\ln\ln k - \frac{1}{2}\ln(2\pi) + \psi_k\left(\log_k\Big(\frac{n}{\log_k n}\Big)\right) + \bigo\left(\frac{\ln^2\ln n}{\ln n}\right)\,,
  \end{aligned}
\end{equation}
where $\psi_k$ is a periodic function with period $1$,
\begin{equation}
  \label{eq:psi-fourier}
  \psi_k(x) = \sum_{j=-\infty}^\infty \alpha_j(k) \E^{2\pi ij x}\,.
\end{equation}
The Fourier coefficients are
\begin{equation}
  \label{eq:alpha-j}
  \alpha_j(k) = \frac{1}{\ln k}\,\Gamma\left(\frac{2\pi i j}{\ln k}\right)\,\zeta\left(1+\frac{2\pi i j}{\ln k}\right) \qquad (j \neq 0)
\end{equation}
and
\begin{equation}
  \label{eq:alpha-0}
  \alpha_0(k) = \frac{1}{\ln k}\left(-\gamma_1-\frac{1}{2}\gamma^2 + \frac{1}{12}\pi^2 + \frac{1}{12} \ln^2 k\right)\,,
\end{equation}
where $\gamma=0.5772156649\ldots$ is the Euler constant and
$\gamma_1=-0.0728158454\ldots$ is the first Stieltjes constant.

The Fourier series for $\psi_k(x)$ converges absolutely and uniformely because
the coefficients $\alpha_j(k)$ decay fast enough:
\begin{equation}
  \label{eq:Gamma-decay}
  |\Gamma(it)| = \bigo(|t|^{-\frac12}\E^{-\frac12\pi|t|})
\end{equation}
and
\begin{equation}
  \label{eq:zeta-decay}
  |\zeta(1+it)| = \bigo(\ln|t|)\,.
\end{equation}
Plugging \eqref{eq:de-bruijn-asymptotics} into \eqref{eq:def-c} 
provides us with
\begin{equation}
  \label{eq:our-c-full}
  \ln c_k(n) = \ln (k \ln k) + \Delta\psi_k(n) + \bigo\left(\frac{(\ln\ln n)^2}{\ln n}\right)\,,
\end{equation}
where 
\begin{equation}
  \label{eq:deltapsi}
  \Delta\psi_k(n) = \psi_k\left(\log_k\Big(\frac{n}{\log_k n}\Big)-\log_k\Big(1+\frac{1}{\log_k n}\Big)\right) 
  - \psi_k\left(\log_k\Big(\frac{n}{\log_k n}\Big)\right)\,.
\end{equation}
Intuitively, $\Delta\psi_k(n)$ should vanish for $n\to\infty$, but to be sure
we need to investigate the Fourier series for $\psi_k$ in more detail. 
In particular, we have
\begin{displaymath}
  \begin{aligned}
    |\Delta\psi_k(n)| &= \left|\sum_{j=-\infty}^\infty \alpha_j(k) \exp\left(2\pi i j \log_k\Big(\frac{n}{\log_k n}\Big)\right)\,\left(\exp\left(-2\pi i j\log_k\Big(1+\frac{1}{\log_k n}\Big)\right)-1\right)\right|\\
   &\leq \sum_{j=-\infty}^\infty |\alpha_j(k)|\,\left|\exp\left(-2\pi i j\log_k\left(1+\frac{1}{\log_k n}\right)\right)-1\right|\\
   &\leq 2\pi \sum_{j=-\infty}^\infty |j\,\alpha_j(k)|\,\left|\log_k\left(1+\frac{1}{\log_k n}\right)\right|\,.
  \end{aligned}
\end{displaymath}
In the last line we have used the inequality
\begin{equation}
  |\E^{ix} - 1| \leq |x| \qquad (x\in\mathbb{R})\,.
\end{equation}
Now because of \eqref{eq:Gamma-decay} and \eqref{eq:zeta-decay} we know that
$\sum_j |j \alpha_j(k)| < \infty$, and hence
\begin{equation}
  \label{eq:deltapsi-decay}
  \Delta\psi_k(n) = \bigo(1/\log_k n)\,.
\end{equation}
This concludes our proof of \eqref{eq:our-c}.

\section{Numerical Evaluation}

The recurrence \eqref{eq:def-lagged-fibo} appears in the analysis of
the Karmarkar-Karp differencing algorithm for number partitioning
\cite{kk:08}. In this context we learned that the convergence
to the asymptotic regime can be extremely slow.
We will see that this is also true when we try to probe the asymptotics of
$c_k(n)$ numerically.

\begin{figure}
  \centering
  \includegraphics[width=\linewidth]{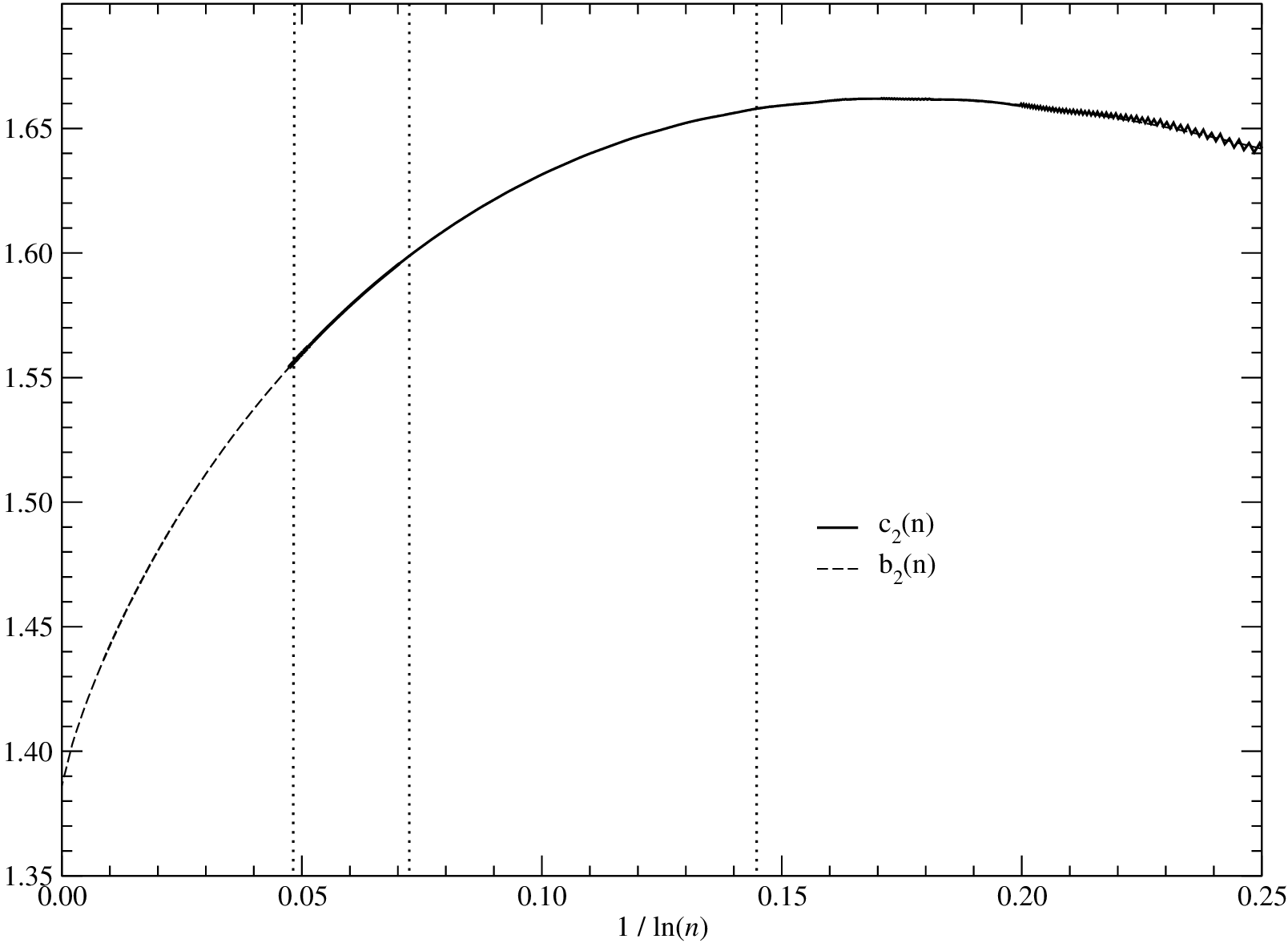}
  \caption{Numerical evaluation of $c_2(n)$ together with $b_2(n)$ from 
    \eqref{eq:def-b}. The difference between $c_2$ and $b_2$ is smaller
    then the linewidth for all $n>100$. The dotted vertical lines
    indicate $n=10^3$, $10^6$ and $10^9$.}
  \label{fig:c}
\end{figure}

To calculate $c_k(n)$, we need $a(kn)$ and $a(n)$, but because of
\begin{equation}
  \label{eq:akn}
  a_k(kn) = a_k(n) + k\sum_{j=0}^{n-1} a_k(j)\,,
\end{equation}
the value of $a_k(n)$ plus the sum of the preceeding terms is
sufficient.  The bottleneck for calculating $a_k(n)$ is memory, not
CPU time, since $n(1-1/k)$ values must be stored to compute $a_k(n)$.
We used the Chinese Remainder Theorem to keep the individual numbers
small and managed to calculate $c_2(n)$ for $n$ up to $3\cdot10^9$ on
a PC with $4$ GByte of memory. As Fig.~\ref{fig:c} shows, even these
data are insufficient to extrapolate to the true asymptotic value.
Numerical calculations that stop at even smaller values of $n$ may
easily misguide an extrapolation to $c_k(\infty)$.

In order to evaluate $c_k(n)$ for much larger values of $n$, we resort
to another asymptotic result. In 1940, Mahler \cite{mahler:40} showed
that
\begin{equation}
  \label{eq:mahler}
  a_k(n) = \E^{\phi_k(n)} \sum_{j=0}^\infty \frac{n^j}{k^{\binom{j}{2}} j!}
  \equiv \E^{\phi_k(n)}\,S_k(n) 
\end{equation}
where $\phi_k(n) = O(1)$. The idea is to replace the numerical
evaluation of $a_k(n)$ by the numerical evaluation of the sum
$S_k(n)$. Note that $S_k(n)$ can be evaluated for very large values of
$n$ using a computer algebra system. A discrepancy in this approach
arises from the unknown function $\phi_k$. Albeit asymptotically
bounded, it can introduce large errors for finite values of $n$.

\begin{figure}
   \centering
   \includegraphics[width=\linewidth]{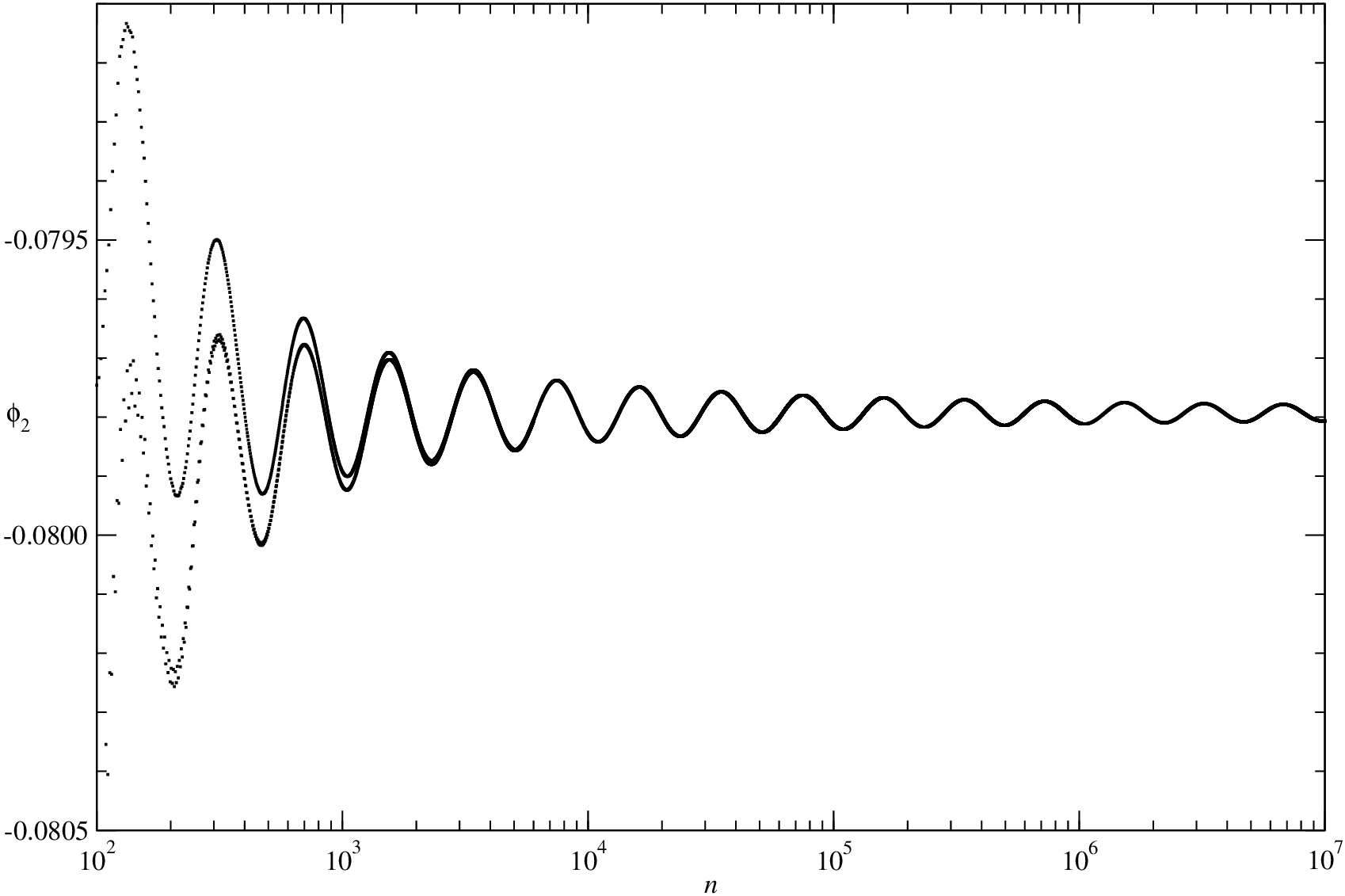}
   \caption{Numerical evaluation of $\phi_2$}
   \label{fig:phi-2}
\end{figure}

\begin{table}
  \centering
   \begin{tabular}{c|c|c}
   $n$ & $c_2(n)$ & $\frac{S_k(kn)}{S_k(n)}\cdot\frac{\ln n}{n}$ \\\hline
   $10^1$ & 1.49668 & 1.50889 \\
   $10^2$ & 1.65470 & 1.65496 \\
   $10^3$ & 1.65791 & 1.65779 \\
   $10^4$ & 1.63881 & 1.63876 \\
   $10^5$ & 1.61782 & 1.61780 \\
   $10^6$ & 1.59883 & 1.59882 \\
   $10^7$ & 1.58237 & 1.58237 \\
   $10^8$ & 1.56822 & 1.56822 \\
   $10^9$ & 1.55600 & 1.55600 \\
 \end{tabular}
  \caption{Exact evaluation of $c_2(n)$ versus evaluation of the series.}
  \label{tab:c-vs-series}
\end{table}

It was already noticed by Fr{\"o}berg \cite{froeberg:77}, that $\phi_k$
oscillates with a small (and decaying) amplitude around a constant
value. We used our extensive data for $a_k(n)$ to look more closely at
\begin{equation}
  \label{eq:def-phi}
  \phi_k(n) \equiv \ln a_k(n) - \ln\sum_{j=0}^\infty \frac{n^j}{k^{\binom{j}{2}} j!}\,. 
\end{equation}
As can be seen from Figure \ref{fig:phi-2}, the amplitude of
$\phi_2$ is smaller than $10^{-4}$ for $n > 10^4$, and it is slowly, but
monotonically decaying. The constant around which $\phi_k$
oscillates will cancel in the ratio
\begin{displaymath}
  \frac{a_k(kn)}{a_k(n)} = \E^{\phi_k(kn)-\phi_k(n)}\,\frac{S_k(kn)}{S_k(n)}
\end{displaymath}
Hence the error is bounded by the small amplitude.  This is confirmed
by the numerical data, see Table~\ref{tab:c-vs-series}. Even for $n=100$, the
error in $c_2(n)$ is only in the fourth decimal. 

This observation tells us that we can use
\begin{equation}
  \label{eq:def-b}
  b_k(n) = \frac{S_k(kn)}{S_k(n)}\,\frac{\ln n}{n} 
\end{equation}
as an excellent approximation to $c_k(n)$. Since $b_k(n)$ can be
evaluated for very large values of $n$, like $n=2^{1000}$ and beyond,
we can use $b_k(n)$ to bridge the gap between the
numerically accessible $c_k(n)$ and $c_k(\infty)$ (Figure~\ref{fig:c}).

\section{Asymptotics Reloaded}

The results of de Bruijn
\eqref{eq:de-bruijn-asymptotics} and Mahler \eqref{eq:mahler} have to
match, i.e., we know that $\phi_k(n) + \ln S_k(n)$ equals the right
hand side of \eqref{eq:de-bruijn-asymptotics}. A saddle-point
expansion of $S_k(n)$ (see \eqref{eq:lnS_Asymp_n} in the Appendix) reveals that
the leading terms of $\ln S_k(n)$ equal the leading terms in
\eqref{eq:de-bruijn-asymptotics}. The remaining terms yield
\begin{equation}
  \label{eq:phi-expansion}
  \phi_k(n) = \Psi_k\left(\log_k\Big(\frac{n}{\log_k n}\Big)\right) +
  \left(1+\frac{\ln\ln k}{\ln k}\right)\ln k - \frac{1}{2}\ln(2\pi) +
  \frac{(\ln k + 2\ln\ln k)^2}{8\ln k} + \bigo\left(\frac{\ln^2\ln n}{\ln n}\right)\,.
\end{equation}

 \begin{figure}
   \centering
   \includegraphics[width=\linewidth]{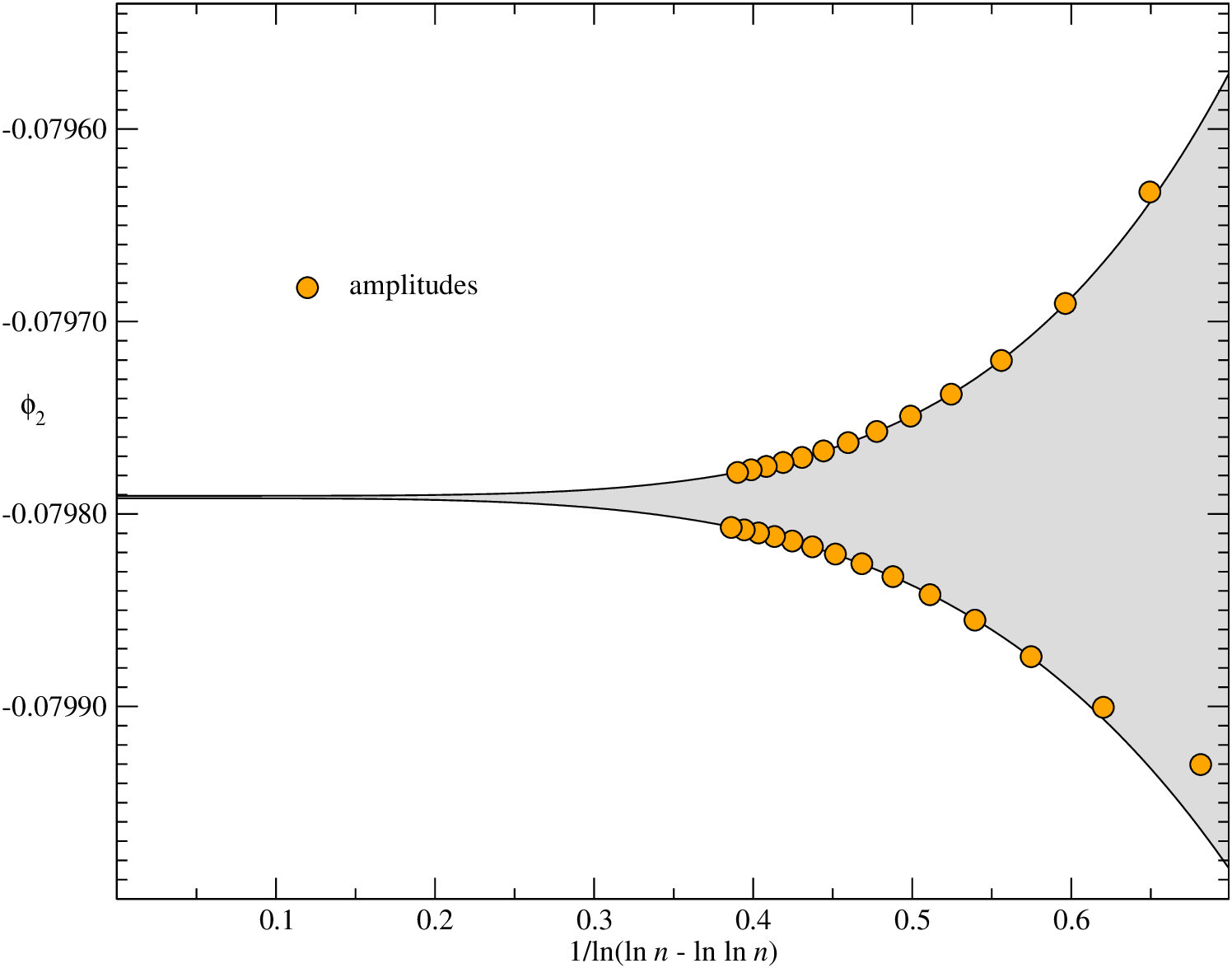}
   \caption{Decay of the amplitude of $\phi_2$. The scale of the abscissa
   is chosen to match $1$ over the period of $\psi_k$. The numerical fits for
   the minima and the maxima are three parameter least square fits
   of $\mu_0 + \mu_1 x^{\mu_2}$.}
   \label{fig:phi-min-max}
 \end{figure}

In particular, we see that asymptotically $\phi_k$ oscillates around a
value
\begin{equation}
  \label{eq:phi-constant}
  \left(1+\frac{\ln\ln k}{\ln k}\right)\ln k - \frac{1}{2}\ln(2\pi) -
  \frac{(\ln k + 2\ln\ln k)^2}{8\ln k} + \alpha_0(k)\,,
\end{equation}
with $\alpha_0(k)$ from \eqref{eq:alpha-0}. For $k=2$, this constant
is $-0.079793025\ldots$ (see Figure~\ref{fig:phi-2}), in perfect
agreement with the numerical results of Fr{\"o}berg
\cite{froeberg:77}.

\begin{table}
  \centering
   \begin{tabular}{c|c|c}
   $j$ & $|\alpha_j(2)|$ & $|\alpha_j(3)|$ \\\hline
   $1$ & $7.36616\cdot10^{-7\phantom{0}}$ &  $1.15010\cdot10^{-4\phantom{0}}$\\
   $2$ & $4.63909\cdot10^{-13}$ & $1.45894\cdot10^{-8\phantom{0}}$ \\
   $3$ & $2.64857\cdot10^{-19}$ & $2.10798\cdot10^{-12}$ \\
   $4$ & $1.29245\cdot10^{-25}$ & $1.45000\cdot10^{-16}$
 \end{tabular}
  \caption{Taylor coefficients \eqref{eq:alpha-j} for $\phi_k$ and $\psi_k$.}
  \label{tab:taylor}
\end{table}

The asymptotic amplitude of $\phi_k$ is very small, as can be seen by
evaluating the coefficients \eqref{eq:alpha-j}, see Table~\ref{tab:taylor}.
Hence we know that the oscillation in Figure \ref{fig:phi-2} will eventually
decay to an amplitude of size $10^{-6}$. We have calculated a few more minima
and maxima of $\phi_2$ to check this decay. Figure~\ref{fig:phi-min-max} shows
the result. The extrapolation of the numerical data gives very accurate result
for the constant $-0.079793025\ldots$ as well as the right order of magnitude
($10^{-6}$) of the remanent amplitude.

\section{Acknowledgments}

We thank Sebastian Mingramm for providing us with the numerical values
for $a_2(n)$ for $n > 10^7$. SB thanks the Institut f{\"u}r Theoretische
Physik at Otto-von-Guericke University in Magdeburg for its
hospitality during the preparation of this manuscript and gratefully
acknowledges support from a Fulbright-Kommision grant and from the 
U.S. National Science Foundation through grant number DMR-0812204.

\section{Appendix}

Here, we evaluate the asymptotic behavior of the sum \begin{eqnarray}
S_{k}(n)=\sum_{j=0}^{\infty}\frac{n^{j}}{j!\, k^{\frac{j(j-1)}{2}}}\qquad(n\to\infty)\label{eq:sum}\end{eqnarray}
 using a saddle-point expansion. Following Ref.~\cite{BO} (pp. 304),
we define $\Phi_{j}=\ln a_{j}$ for the summand $a_{j}$. The finite-difference
condition $D\Phi_{j}=\Phi_{j}-\Phi_{j-1}=0$ determines the maxima,
i.~e. we need to find the $j_{0}$-term(s) of the sum with $a_{j_{0}}/a_{j_{0}-1}\sim1$.
Applied to Eq.~(\ref{eq:sum}), we obtain \begin{eqnarray}
\eta\left[j_{0}+\Delta_{j}-1+\log_{k}\left(j_{0}+\Delta_{J}\right)\right] & \sim & 1\label{eq:saddlep}\end{eqnarray}
 where we use the abbreviations \begin{eqnarray}
\eta=\frac{1}{\log_{k}n} & = & \frac{\ln k}{\ln n}\label{eq:def}\end{eqnarray}
 and a non-integer offset $0\leq\Delta_{j}\leq1$ on the integer location
$j_{0}$ of the saddlepoint, which we will need to attain the continuum
limit for this $n$-dependent ({}``moving'') saddle point \cite{BO}.
For $n\to\infty$, $\eta\to0$ and in that limit we find from Eq.
(\ref{eq:saddlep}) for the saddle-point location by peeling off layer-by-layer:
\begin{eqnarray}
j_{0} & = &
\frac{1}{\eta}+\log_{k}\eta+1-\Delta_{j}-\eta\frac{1+\log_{k}\eta}{\ln
  k}+\eta^{2}\frac{1+\log_{k}\eta}{2\ln  k}
\left[\frac{2}{\ln
    k}+1+\log_{k}\eta\right]\nonumber\\
&&\quad -\eta^3\frac{1+\log_k\eta}{6\ln k}
\left[\frac{6}{\ln^2k}+\frac{9}{\ln k}+2+\frac{\ln\eta}{\ln
    k}\left(4+\frac{9}{\ln k}\right)+2\frac{\ln^2\eta}{\ln^2k}\right]
+O\left(\eta^{4}\ln^{4}\eta\right).
\label{eq:sp}
\end{eqnarray}
 As $j_{0}\gg1$ for $\eta\to0$, we can expand $\Phi_{j}$ for large
arguments: \begin{eqnarray}
\Phi_{j} & = & j\ln\left(n\right)-\ln\left(j!\right)-\frac{j\left(j-1\right)}{2}\ln\left(k\right),\label{eq:asymp_phi}\\
 & = & \frac{j}{\eta}\ln\left(k\right)-\frac{j\left(j-1\right)}{2}\ln\left(k\right)+\left(j+\frac{1}{2}\right)\ln\left(j\right)+j-\frac{1}{2}\ln\left(2\pi\right)+\frac{1}{12j}+O\left(\frac{1}{j^{3}}\right),\nonumber \end{eqnarray}
 where we have used the Stirling expansion for the factorial to the
necessary order.

At the (unique) maximum of $\Phi_{j}$ we set $j\sim j_{0}+t$ and
expand here only to quadratic (Gaussian) order in $t$:%
\footnote{Higher orders in $t$ are irrelevant here for the order in $\eta$
considered.%
} \begin{eqnarray}
\Phi_{j_{0}+t} & = & \frac{j_{0}}{\eta}\ln\left(k\right)-\frac{j_{0}\left(j_{0}-1\right)}{2}\ln\left(k\right)+\left(j_{0}+\frac{1}{2}\right)\ln\left(j_{0}\right)+j_{0}-\frac{1}{2}\ln\left(2\pi\right)+\frac{1}{12j_{0}}\nonumber \\
 &  & \quad+t\left[\left(\frac{1}{\eta}-j_{0}+\frac{1}{2}\right)\ln\left(k\right)-\frac{\ln\left(j_{0}\right)}{2}-\frac{1}{2j_{0}}+\frac{1}{12j_{0}^{2}}\right]-t^{2}\left[\frac{\ln\left(k\right)}{2}+\frac{1}{2j_{0}}-\frac{1}{4j_{0}^{2}}\right]+O\left(\frac{t^{3}}{j_{0}^{2}}\right).\label{eq:t-expansion}\end{eqnarray}
 Note that the linear term in $t$ only vanishes (indicating a symmetric
maximum) after we insert the moving saddle-point in Eq. (\ref{eq:sp})
and $\Delta_{j}$ is fixed: \begin{eqnarray}
\Phi_{j_{0}+t} & = & \frac{\ln\left(k\right)}{2\eta^{2}}+\frac{2\ln\left(\eta\right)+\ln\left(k\right)+2}{2\eta}+\left[\frac{\ln^{2}\left(\eta\right)}{2\ln\left(k\right)}+\ln\left(\eta\right)-\frac{1}{2}\ln\left(2\pi\right)+\frac{\Delta_{j}-\Delta_{j}^{2}}{2}\ln\left(k\right)\right]\nonumber \\
 &  & \quad-\eta\left[\frac{\ln^{2}\left(\eta\right)}{2\ln^{2}\left(k\right)}+\frac{\ln\left(\eta\right)}{\ln\left(k\right)}+\frac{7}{12}-\frac{\Delta_{j}-\Delta_{j}^{2}}{2}\right]\nonumber \\
 &  & \quad+\frac{\eta^{2}}{12}\left[\left(1+\frac{\ln\eta}{\ln k}\right)\left(6\frac{1+\frac{\ln\eta}{\ln k}}{\ln k}+3+4\frac{\ln\eta}{\ln k}+2\frac{\ln^{2}\eta}{\ln^{2}k}\right)+\left(\Delta_{j}-\Delta_{j}^{2}\right)\left(2\Delta_{j}-7-6\frac{\ln\eta}{\ln k}\right)\right]\\
 &  & \quad+O\left(\eta^{3}\ln^{4}\eta\right)\\
 &  & \quad+\left[\frac{1}{2}\left(\ln\left(k\right)+\eta\right)\left(1-2\Delta_{j}\right)+\eta^{2}\left(\frac{7}{12}-\frac{3}{2}\Delta_{j}+\frac{1}{2}\Delta_{j}^{2}+\frac{1-2\Delta_{j}}{2}\,\frac{\ln\eta}{\ln k}\right)\right]t\label{eq:full-expansion}\\
 &  & \quad-\frac{1}{2}\left[\left(\ln\left(k\right)+\eta\right)\left(1-2\Delta_{j}\right)+\eta^{2}\left(\Delta_{j}-\frac{3}{2}-\frac{\ln\eta}{\ln k}\right)\right]t^{2}+O\left(\eta^{2}t^{3}\right).\end{eqnarray}
Terms of orders such as $\eta^{2}t^{3}$ will not contribute at order
$\eta^{2}$ as they are at leading order asymmetric in $t$ in the
ensuing Gaussian integration. To that effect, we symmetrize the saddle
point to order $\eta^{2}$ with the choice of \begin{eqnarray}
\Delta_{j} & = & \frac{1}{2}+\frac{\eta^{2}}{24\ln k}+O\left(\eta^{3}\ln\eta\right),\label{eq:Delta_choice}\end{eqnarray}
which also impacts constant or smaller terms in $\eta$ in Eq. (\ref{eq:full-expansion}).
The Gaussian integration then yields \begin{eqnarray}
S_{k}\left(n\right) & \sim & \sum_{t=-\epsilon j_{0}}^{\epsilon j_{0}}e^{\Phi_{j_{0}+t}},\qquad\left(\eta\ll\epsilon\ll1\right)\label{eq:S-asymp}\\
 & = & e^{\Phi_{j_{0}}}\int_{-\infty}^{\infty}dt\,\exp\left\{ -\frac{1}{2}\left[\ln\left(k\right)+\eta-\eta^{2}\left(1+\frac{\ln\eta}{\ln k}\right)\right]t^{2}+O\left(\eta^{2}t^{3}\right)\right\} ,\nonumber \\
 & = & e^{\Phi_{j_{0}}}\int_{-\infty}^{\infty}dt\,\left[1+O\left(\eta^{2}t^{3}\right)\right]\exp\left\{ -\frac{1}{2}\left[\ln\left(k\right)+\eta-\eta^{2}\left(1+\frac{\ln\eta}{\ln k}\right)\right]t^{2}\right\} ,\\
 & = & \exp\left\{ \Phi_{j_{0}}+\frac{1}{2}\ln\left(2\pi\right)-\frac{1}{2}\ln\left[\ln\left(k\right)+\eta-\eta^{2}\left(1+\frac{\ln\eta}{\ln k}\right)\right]+O\left(\eta^{3}\ln^{2}\eta\right)\right\} .\nonumber \end{eqnarray}
 Note that the $\ln\left(2\pi\right)$-terms cancel. Hence, we finally
obtain\begin{eqnarray}
\ln\left[S_{k}\left(n\right)\right] & = & \frac{\ln k}{2\eta^{2}}+\frac{2\ln\eta+\ln k+2}{2\eta}+\left[\frac{\ln^{2}\eta}{2\ln k}+\ln\eta-\frac{1}{2}\ln\ln k+\frac{1}{8}\ln k\right]\nonumber \\
 &  & \quad-\eta\left[\frac{\ln^{2}\eta}{2\ln^{2}k}+\frac{\ln\eta}{\ln k}+\frac{11}{24}+\frac{1}{2\ln k}\right]\nonumber \\
 &  & \quad+\eta^{2}\left[\frac{\ln^{3}\eta}{6\ln^{3}k}+\left(1+\frac{1}{\ln k}\right)\frac{\ln^{2}\eta}{2\ln^{2}k}+\left(\frac{11}{24}+\frac{3}{2\ln k}\right)\frac{\ln\eta}{\ln k}+\left(\frac{1}{8}+\frac{1}{\ln k}+\frac{1}{4\ln^{2}k}\right)\right]\nonumber \\
 &  & \quad+O\left(\eta^{3}\ln^{4}\eta\right),\label{eq:lnS_asymp}\end{eqnarray}
or in terms of powers of $\ln n$ directly:\begin{eqnarray}
\ln\left[S_{k}\left(n\right)\right] & = & \frac{\ln^{2}n}{2\ln k}+\ln n\left[-\frac{\ln\ln n}{\ln k}+\frac{\ln\ln k}{\ln k}+\frac{1}{2}+\frac{1}{\ln k}\right]+\frac{\ln^{2}\ln n}{2\ln k}-\ln\ln n\left[\frac{\ln\ln k}{\ln k}+1\right]+\frac{\left[\ln k+2\ln\ln k\right]^{2}}{8\ln k}\nonumber \\
 &  & \quad-\frac{\ln k}{\ln n}\left[\frac{\ln^{2}\ln n}{2\ln^{2}k}-\frac{\ln\ln n}{\ln k}\left(1+\frac{\ln\ln k}{\ln k}\right)+\frac{11}{24}+\frac{1}{2\ln k}+\frac{\ln\ln k}{\ln k}+\frac{\ln^{2}\ln k}{2\ln^{2}k}\right]\nonumber \\
 &  & \quad-\frac{\ln^{2}k}{\ln^{2}n}\left[\frac{\ln^{3}\ln n}{6\ln^{3}k}-\frac{\ln^{2}\ln n}{2\ln^{2}k}\left(1+\frac{1}{\ln k}+\frac{\ln\ln k}{\ln k}\right)\right.\nonumber \\
 &  & \qquad+\frac{\ln\ln n}{\ln k}\left(\frac{11}{24}+\frac{3}{2\ln k}+\frac{\ln\ln k}{\ln k}+\frac{\ln\ln k}{\ln^{2}k}+\frac{\ln^{2}\ln k}{2\ln^{2}k}\right)\nonumber \\
 &  & \qquad\left.-\frac{1}{8}-\frac{1}{\ln k}-\frac{1}{4\ln^{2}k}-\frac{11\ln\ln k}{24\ln k}-\frac{3\ln\ln k}{2\ln^{2}k}-\frac{\ln^{2}\ln k}{2\ln^{2}k}-\frac{\ln^{2}\ln k}{2\ln^{3}k}-\frac{\ln^{3}\ln k}{6\ln^{3}k}\right]\nonumber \\
 &  & \quad+O\left(\frac{\ln^{4}\ln n}{\ln^{3}n}\right).\label{eq:lnS_Asymp_n}\end{eqnarray}

\begin{figure}
\includegraphics[width=\linewidth]{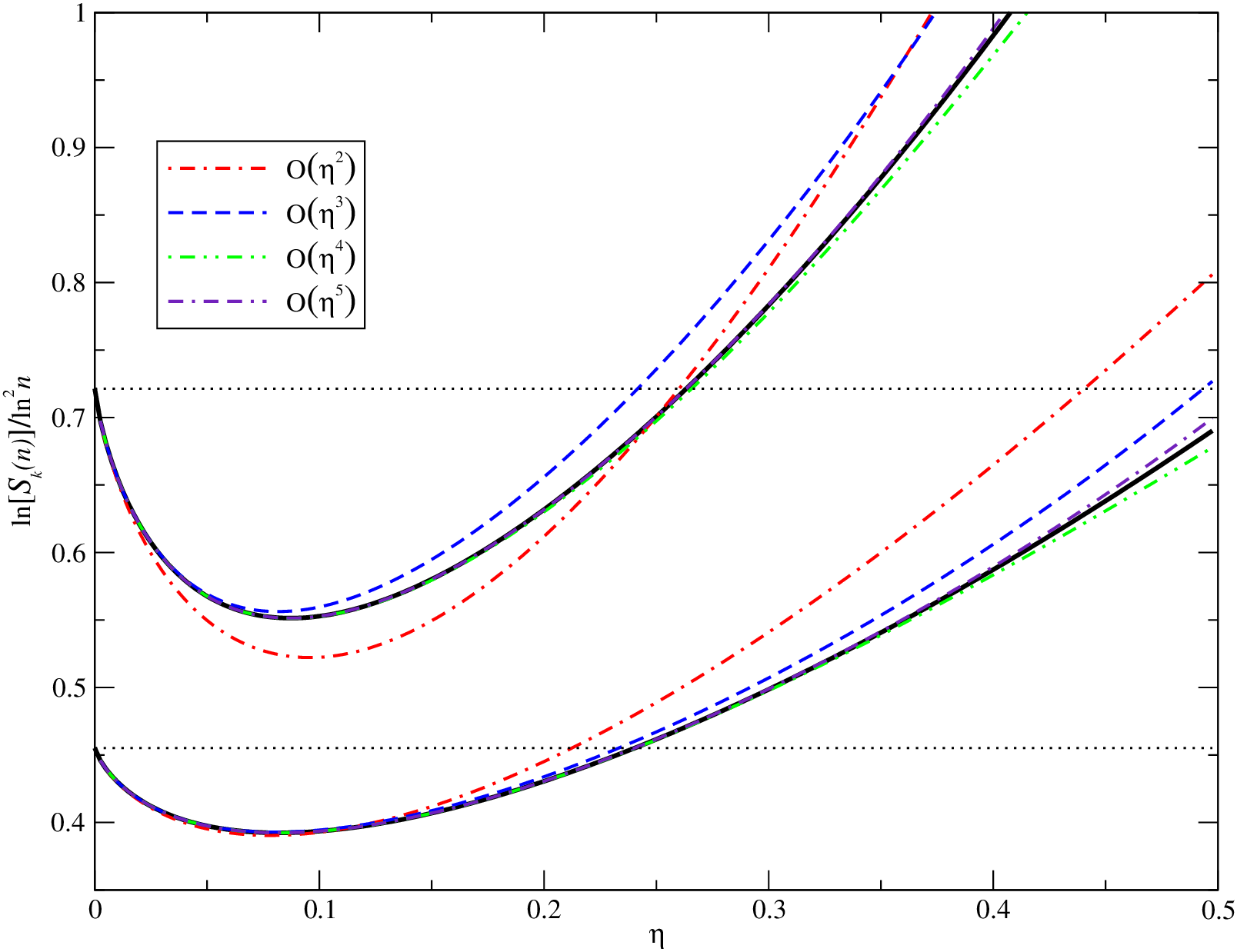}\caption{\label{fig:Plot-of-lnS}Plot of $\ln S_{k}(n)/\ln^{2}n$ vs $\eta=\ln k/\ln n$
for $k=2$ (top set) and 3 (bottom set). The dotted horizontal lines
specify the asymptotic limits, $\frac{1}{2\ln2}=0.7213\ldots$ for
$k=2$ and $\frac{1}{2\ln3}=0.4551\ldots$ for $k=3$. The thick black
line is obtained from the numerically exact evaluation of Eq. (\ref{eq:sum}),
and the shaded, dashed lines correspond to the asymptotic expression
in Eq. (\ref{eq:lnS_asymp}), evaluated to the indicated order in
$\eta$.}

\end{figure}

In Fig. \ref{fig:Plot-of-lnS} we plot a sequence of approximants
to the numerically exact evaluation of the sum in Eq. (\ref{eq:sum}),
which prove to approximate with an error of the indicated order. 

\bibliography{mertens,fibo,math}

\bigskip
\hrule
\bigskip

\noindent 2000 {\it Mathematics Subject Classification}: 
Primary 11B39; Secondary 41A60, 11P81

\noindent \emph{Keywords: } 
partitions, Fibonacci sequence, linear
recurring sequence, asymptotics, Mahlerian sequence

\bigskip
\hrule
\bigskip

\noindent (Concerned with sequence
\seqnum{A033485}, 
\seqnum{A000123}, 
\seqnum{A005704}, 
\seqnum{A005705}, and 
\seqnum{A005706}.) 

\bigskip
\hrule
\bigskip

\end{document}